\documentclass[11pt,a4paper]{amsart}
\usepackage{amsmath}
\usepackage{amsthm}
\usepackage{amsfonts}
\usepackage{amsfonts,amsmath,amstext,amsbsy,euscript,amssymb, graphics}
\usepackage{rotating}
\usepackage{enumerate}
\font\amsx=msam10 \font\amsy=msbm10
\def\IZ{\hbox{\amsy\char'132}}
\def\R{\mathbb{R}}
\def\N{\mathbb{N}}
\def\NN{\mathbb{N}_{0}}
\def\M{\mathcal{M}}
\def\P{\mathcal{P}}
\def\n{\mathcal{N}}
\def\B{\mathcal{B}}
\def\T{\mathcal{T}}
\def\W{\text{W}}
\newtheorem{theorem}{Theorem}
\newtheorem{lemma}{Lemma}
\newtheorem{cor}{Corollary}
\newtheorem{prop}{Proposition}
\newtheorem{conj}{Conjecture}

\newtheorem{rem}{Remark}

\newtheorem{ques}{Question}
\newenvironment{pf}{\medskip\noindent{\bf Proof}.}{}
\def\bqed{\quad\hbox{\amsx\char'004}}

\title{The $\star$-operator and Invariant Subtraction Games}
\author{Urban Larsson}
\begin{document}
\keywords{Beatty sequence, Complementary sequences, Dual game, 
Fibonacci sequence, Impartial game, Invariant subtraction games, Involution, 
Nim, Permutation, $\star$-operator, Wythoff Nim, Zeckendorf numeration}
\begin{abstract}
We study 2-player impartial games, so 
called \emph{invariant subtraction games}, of the type, given a set 
of allowed moves the players take turn in moving one single piece on a 
large Chess board towards the position $\boldsymbol 0$. Here, invariance 
means that each allowed move is available inside the whole board. 
Then we define 
a new game, $\star$ of the old game, by taking the 
$P$-positions, except $\boldsymbol 0$, as moves in 
the new game. One such game is $\W^\star=$ (Wythoff Nim)$^\star$, where the 
moves are defined by complementary Beatty sequences with irrational moduli. 
Here we give a polynomial time algorithm for infinitely many $P$-positions 
of $\W^\star$. A repeated application of 
$\star$ turns out to give especially nice properties for 
a certain subfamily of the invariant subtraction games, 
the \emph{permutation games}, which we introduce here.
We also introduce the family of \emph{ornament games}, 
whose $P$-positions define 
complementary Beatty sequences with rational moduli---hence related 
to A.S Fraenkel's `variant' Rat- and Mouse games---and give closed forms 
for the moves of such games. 
We also prove that ($k$-pile Nim)$^{\star\star}$ = $k$-pile Nim. 

\end{abstract}
\maketitle
\section{Introduction and terminology}
This paper is a sequel to \cite{LHF}. We begin by recapitulating 
some terminology. A 2-player \emph{impartial} \cite{BCG1982} 
game is a combinatorial 
game where, independent of whose turn it is, the options are the same. 
Here (as in \cite{LHF}) we study so called \emph{invariant subtraction games}\footnote{The subtraction games in \cite{BCG1982} are special cases 
of the ones studied here.}, impartial `board games', mostly played 
on the board $\B = \NN\times \NN$ (except in the last section where 
it is $\NN^k$, $k\in \N )$. Given a set of `invariant moves', 
denoted by $\M (G)$, the two players take turn in moving a single piece towards 
the position $\boldsymbol 0 = (0, 0)$. The player who moves there 
wins. A `move' is represented by an ordered pair of non-negative integers, 
say $(i, j)\ne \boldsymbol 0$. In practice, this ordered pair is 
subtracted from the piece's current position, say $(x,y)$. The 
resulting position of this move, provided it is allowed, 
is $$(x,y)\ominus (i,j)\succeq \boldsymbol 0.$$
Here \emph{invariance} (of the the set of moves) means that  
each allowed move is playable from any position of $\B$, 
provided that the piece remains on the board.
 
Nim \cite{B1902} is a classical impartial game played 
on a finite number of piles each containing a finite number of tokens. 
The players take turn in removing tokens from precisely one of the piles, 
at least one token and at most the whole pile.
It may be regarded as an invariant subtraction game 
with, if played on two piles, 
$$\M (\text{2-pile Nim}) = \{\{0,x\} \mid x\in \N\}.$$
(We use the `symmetric notation' $\{x, y\}$ whenever the ordered 
pairs $(x,y)$ and $(y,x)$ are considered the same.) Another classical 
example of an invariant subtraction game is Wythoff Nim \cite{W1907}, 
here denoted by W. The players take turn in moving a Queen of Chess on a 
large Chess board towards the lower-left corner. In our notation, the moves are 
$$\M(\text{W}) = \M(\text{2-pile Nim})\cup \{(x, x)\mid x\in \N\}.$$ 

As many impartial games, invariant subtraction games have no draw 
(cyclic) moves and hence the positions 
are either $P$ (the previous player wins) or $N$ (the next player wins). 
Given a game $G$, the sets of all $P$-positions and all
$N$-positions is denoted by $\P(G)$ and $\n(G)$ respectively. 
We denote the set of \emph{terminal} positions by 
$T(G)\subset \P(G)$. It contains all positions with empty sets 
of options.

\subsection{Non-zero $P$-positions as moves}\label{sec:1.1}
We may now define the $\star$-operator, introduced in \cite{LHF}.
Suppose that $G$ is an invariant subtraction game. Then $G^\star$ is 
the game defined by $$\M(G^\star) = \P(G)\setminus \{\boldsymbol 0\}.$$
We let $G^k$ denote the resulting 
game of $k$ recursive applications of $\star$, so that for example $G^0=G$ and
$G^3 = ((G^\star)^\star)^\star$. For the special case of $k=2$ we prefer to 
write $G^{\star\star}$. As in \cite{LHF}, if 
$G = G^{\star\star}$ we say that $G^\star$ is the \emph{dual} of $G$. 

\subsection{Games defined via complementary Beatty sequences}\label{sec:1.2}
Two sequences of positive integers $(a_i)$ and $(b_i)$ 
are \emph{complementary} 
if $\{a_i\}\cup\{b_i\} = \N$ and $\{a_i\}\cap\{b_i\} = \emptyset$.
A \emph{Beatty sequence} is a sequence of 
the form $(\lfloor \alpha n + \gamma\rfloor)$, where $\alpha, \gamma \in \R$ 
and where $n$ ranges over $\N$. Suppose we have a pair of Beatty sequences, 
say 
\begin{align}\label{beatty1}
(\lfloor \alpha n + \delta\rfloor ) 
\text{ and } (\lfloor \beta n + \gamma\rfloor ). 
\end{align}
Necessary and sufficient 
conditions on their respective moduli and offsets for them to 
be complementary are given in \cite{F1969, OB2003}.

\begin{figure}[ht!]
\centering
\includegraphics[width=0.9\textwidth]{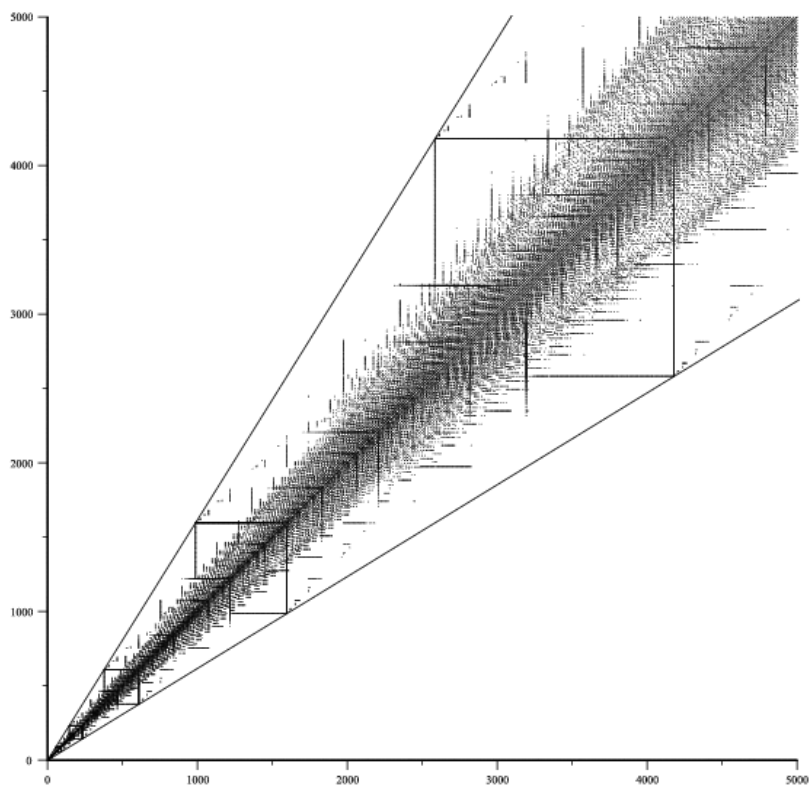}
\caption{The $P$-positions of (Wythoff Nim)$^\star$ with coordinates 
less than 5000 together with the lines through the origin 
with slopes $\phi$ and $\phi^{-1}$ respectively. 
(Remark: There are no $P$-positions `on' these lines.) See 
also Table \ref{table:1}. In Theorem \ref{theorem:2} we prove the existence 
of infinitely many `$\log$-periodic' $P$-positions.
}\label{figure:1}
\end{figure}

\begin{figure}[ht!]
\centering
\includegraphics[width=0.7\textwidth]{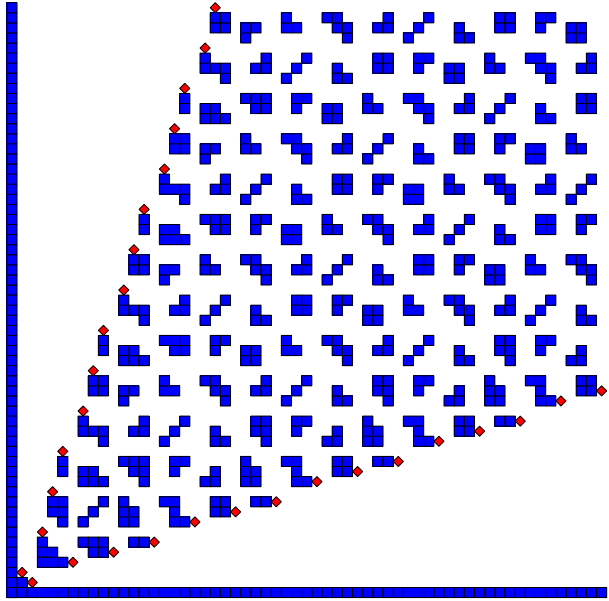}
\caption{ For this game, the red diamonds represent moves given 
by the pairs of complementary Beatty sequences 
$(\lfloor \frac{4n+1}{3}\rfloor)$ 
and $(4n-2)$, where $n$ runs over the positive integers. By this we mean that, 
given a position $(x,y)\in \B$, each legal option is of the form  
 $(x,y)\ominus \{\lfloor \frac{4n+1}{3}\rfloor, 4n-2\}\succeq \boldsymbol 0$.
The apparent 
`periodicity' of the $P$-positions (which stands in bright contrast to  
the graph of the $P$-positions of $\W^\star$ 
see Figure \ref{figure:1}) motivates a `periodicity'-conjecture 
on games defined by complementary Beatty sequences (Conjecture 
\ref{conj:4}). See also Figure \ref{figure:2}, \ref{figure:3} 
and \ref{figure:4}. By \cite[Main Theorem]{LHF}, 
the dual definition is that the blue squares, 
except $\boldsymbol 0$, represent the first few moves of a game 
where the red diamonds represent 
the first few non-zero $P$-positions. Indeed this gives 
one of the four Class 4 ornament games defined in Section 5.}\label{figure:5}
\end{figure}

It is well-known that the set of $P$-positions of Wythoff Nim may be 
defined via complementary Beatty sequences with irrational moduli, namely
$\P(\W)=\{\{\lfloor \frac{\sqrt{5}-1}{2}n\rfloor, 
\lfloor \frac{\sqrt{5}+1}{2}n\rfloor\}\mid n\in \NN\}$.  
We give a polynomial time algorithm for infinitely many $P$-positions of 
the dual $\W^\star$ \cite[Main Theorem]{LHF} of Wythoff Nim (which 
corresponds to infinitely many moves of $\W^{\star\star}$) see also 
Figure \ref{figure:1}. 
We give a closed formula for the set of moves of the 
invariant subtraction game, `the Mouse trap' \cite{LHF}.  
Here the set of $P$-positions 
$\{\{\lfloor \frac{3n}{2}\rfloor, 3n-1\}\mid n\in \N\}$ is defined via 
complementary Beatty sequences with rational moduli. 
(Thus, this game has the same $P$-positions as the `variant' Mouse game 
introduced in \cite{F2010}). We present some more general results on the 
family of all invariant subtraction games for which the sets of $P$-positions 
consists of Beatty sequences with rational moduli, here 
we introduce the notion of \emph{ornament games} 
(e.g. Figure \ref{figure:5}). 

We count the number of such games contained in certain classes and thereby 
demonstrate that, in total, there are only countably many ornament games. 
(In contrast, as noted already in \cite{LHF}, there are uncountably 
many invariant subtraction games with 
sets of $P$-positions defined via irrational Beatty sequences.)
Then, we state a conjecture on invariant subtraction games defined 
via complementary Beatty sequences saying the
the set of $P$-positions is `periodic'
if and only if the moduli of the respective Beatty sequences is rational.

\subsection{Permutation and involution games}\label{sec:1.3}
Let us here introduce the notion of 
a \emph{permutation} game. This is an invariant subtraction game, 
where each row 
and column of $\N\times \N$ contain precisely one move, but, where 
both row 0 and column 0 are void of moves. 
We say that a set, say $S\subset \B$, is 
symmetric if $(x, y)\in S$ if and only if $(y, x)\in S$. 
An \emph{involution} game is a permutation game where the set of moves 
is symmetric. Our main results on permutation games are the following. 
Let $G$ denote a permutation game. Then 
$G^{\star\star}$ is also, so the permutation games are closed under the 
operation $\star\star$. In fact, even more is true. The sequence 
$(G^{2k})_{k\in \N}$ `converges' and, by the closure property, 
the resulting game is a permutation game. Similar results hold for 
involution games.\footnote{These 
results were posed as questions on a seminar I gave in the Spring 2010, 
first at CANT 2010 and then at Dalhousie University. \cite{seminar}}

\subsection{Nim and its dual}\label{sec:1.4}
The $P$-positions of Nim on $k$ piles can be taken as moves in a new game, 
Nim$^\star$. We prove that the non-zero $P$-positions of Nim$^\star$ are 
the moves of Nim. A side effect of this result is a winning strategy 
of Nim without the mention of `Nim sum'. 

\subsection{Exposition} In Section 2 we prove some very basic results for 
invariant 2-pile subtraction games on complementary sequences 
of positive integers.
In Section 3 we discuss the winning strategy of $\W^\star$. In Section 4 
we study permutation games. In Section 5 we 
study the invariant subtraction game `the Mouse trap', including a relative 
to this game with a so-called `Muller twist', complementary 
Beatty sequences with rational moduli and the family of ornament games.
In Section 6 we study Nim$^\star$.

\section{Invariant subtraction games defined by complementary sequences}
We begin with a very basic result concerning 
invariant subtraction games for which the 
set of moves is defined via complementary sequences of positive integers. 

\begin{theorem}\label{theorem:1}
Let $a=(a_i)_{i\in \N}$ and $b=(b_i)_{i\in \N}$ denote complementary sequences of 
positive integers, $a$ increasing, and for all $i$, $a_i < b_i$. 
Define $G$ by $\M(G) = \{\{a_i, b_i\}\mid i\in \N \}$. 
Then, 
\begin{enumerate}[(i)]
\item $(x, y)\in \P(G)\setminus \T(G)$ implies that 
there is an $i\in \N$ such that $x = a_i$ or $y = a_i$.
\item if $b$ is increasing, then $(x, y)\in \P(G)\setminus \T(G)$ 
implies that there are $i, j\in \N$ such that $x = a_i$ and $y = a_j$.
\item if $b_i/a_i$ is bounded by some constant, say $C\in \R$, then 
$(x, y)\in \mathcal{P}(G)\setminus \T(G)$ 
(with $x\le y$) implies that $y/x\le C$.
\end{enumerate}
\end{theorem}

\begin{pf}
By definition of $G$, all positions of the form $\{0,x\}$, $x\in \NN$, 
are terminal and hence $P$.\\

\noindent Case (i). Suppose $i, j \in \N$ with $i\le j$. 
Then $ (b_i, b_j)\ominus (b_i, a_i)$ is terminal, 
hence $\{b_i, b_j\}$ is $N$. The claim follows by complementarity 
of $a$ and $b$.\\\\ 
\noindent Case (ii). 
If the claim does not hold, then there is an $i\in \N$ is such that
\begin{enumerate}[(a)]
\item $y = b_i,$ or
\item $x = b_i,$
\end{enumerate}
But the game is symmetric so it suffices to investigate Case (a).
Suppose, in addition, that $x \ge a_i$. 
Then the option $(x, b_i)\ominus (a_i, b_i)$ 
is terminal. Otherwise, by complementarity, 
there is a $j < i$ such that either $a_j = x$ or $b_j = x$. 
If $b_j = x$, then $b_i > a_j$ so that the option 
$(b_j, b_i)\ominus (b_j, a_j)$ is terminal.
Suppose rather that $x = a_j$. Then, since $b$ is increasing, the option 
$(a_j, b_i)\ominus (a_j, b_j)$ is legal and hence terminal. 
By symmetry we may conclude that no position of the form $\{x, b_i\}$ can 
be $P$.\\\\ 
\noindent Case (iii). Suppose that $y/x > C$. Then, by complementarity, there 
is an $i$ such that either $b_i = x$ or $a_i = x$. In the first case, 
$(x, y)\ominus (b_i, a_i)$ is the desired terminal option and in the second 
case, since $y\ge b_i$, $(x, y)\ominus (a_i, b_i)$ is.
We are done. \hfill \bqed
\end{pf}\medskip

\section{A polynomial time algorithm for infinitely 
many $P$-positions of $\W^{\star}$}

Let $\phi := \frac{1+\sqrt{5}}{2}$ denote the golden ratio 
and, for all $n\in \NN$, define $$A_n := \lfloor \phi n \rfloor$$ 
and $$B_n:= A_n + n.$$ We keep this notation for the rest of this section. 
Then $\P(\W) = \{\{A_i, B_i\}\mid i\in \NN\}$ \cite{W1907}. 
Thus we have a polynomial time algorithm in $\log n$ 
for Wythoff Nim's decision problem: 
Determine whether a given pair of natural numbers 
represents a $P$-position. 

Let $F_0 = F_1 = 1, F_n = F_{n-1} + F_{n-2},$ $(n\ge 2)$ 
denote the sequence of Fibonacci numbers. The main result of this section is.

\begin{theorem}\label{theorem:2}
For all $n$, provided both coordinates are positive,
the following positions (and its symmetric counterparts) 
belong to $\P(\W^\star)$ 
\begin{enumerate}[(i)]
\item $(F_{2n} - 1, F_{2n} - 1)$,
\item $(F_{2n-1}, F_{2n} - 1)$,
\item $(F_{2n-1}, F_{2n} - 4)$,
\item $(F_{2n-1}, F_{2n} - 9)$,
\item $(F_{2n-1} + 1, F_{2n} - 1)$,
\item $(F_{2n-1} + 3, F_{2n} - 1)$,
\item $(F_{2n-1} + 4, F_{2n} - 1)$ and
\item $(F_{2n-1} + 6, F_{2n} - 1)$.
\end{enumerate}
\end{theorem}

In this section we make frequent use of the Fibonacci numeration system. 
Namely, each non-negative 
integer can be represented as a sum of distinct Fibonacci numbers. 
Hence we may code any non-negative integer by some binary string 
$$\alpha_n\alpha_{n-1}\ldots \alpha_1 := \sum_{i = 1}^{n} \alpha_i F_i,$$
for some $n\in \N$ and where, for all $i$, $\alpha_i\in\{0,1\}$.

The Zeckendorf numeration system is the unique Fibonacci numeration, 
where the binary string contains no two consecutive ones.

However, certain properties of a number 
do not depend on which Fibonacci system of numeration we have used.

\begin{prop}\label{prop:1}
Let $X\in \N$. Then a Fibonacci coding of $X$, $X_{fib}$, ends 
in an even number of 0s if and only if its Zeckendorf coding, $X_{zeck}$, does. 
\end{prop}

\begin{pf} Search the digits of $X_{fib}$ from left to right. 
Whenever two consecutive 1s are detected exchange ``11'' for ``100'' (where 
the least 0 has the same position as the previous least 1). 
Repeat this step until no more ``11''s are detected. Then the parity 
of the number of rightmost 0s in $X_{fib}$ is the same as in 
the `output', $X_{zeck}$.
\hfill \bqed
\end{pf}
\medskip

\begin{lemma}[\cite{F1982}] \label{Fr82}
Let $A_n$ and $B_n$ be defined as above. 
Then, in Fibonacci coding, $A_n$ ends in an even number of 0s 
and $B_n = A_n0$. 
\end{lemma}

Combining this result with Theorem \ref{theorem:1} (ii) we obtain 
the following nice property for the strategy of $\W^\star$. This result 
was first proved by A. S. Fraenkel.

\begin{cor}[A.S.Fraenkel]\label{cor:1}
Suppose $(x, y)$ represents a $P$-position of $\W^\star$. Then, in 
Fibonacci coding, both $x$ and $y$ end in an even number of 0s.
\end{cor}

In itself, this result does not reduce the complexity 
of the decision problem for $\W^\star$ to polynomial 
time\footnote{For small $n$ ($\le 50000$), 
it does give a considerable improvement in computing capacity. In fact, 
it seems that then the bounds on the memory (storage of $P$-positions) 
sets the limit of computation rather than the processing power.}
in $\log n$. However, it is clear that it 
characterizes a substantial fraction of the $N$-positions 
in polynomial time. In fact, as we will see, 
this constitutes one of the primary tools 
for proving polynomial complexity of \emph{certain} $P$-positions. 

The next two results concerns arithmetical properties 
of numbers of the form $F_{2n}-1$ and $F_{2n-1}$ respectively.
In the `Fibonacci coding' of a number, we let the symbol $0^t$ 
denote a repetition of $t$ consecutive 0s or, for that matter, we let 
$x^t$ denote a consecutive repetition of $t$ $x$:s (for example 
$F_{12} + F_6 + F_4 + F_2 =100000101010 = 10^5(10)^3$).

\begin{lemma}\label{lemma:1} Let $n\in \N$ and let $X\in \N$ be such that 
$F_{2n} - X > 0$. In Fibonacci coding, put 
$$\xi := F_{2n} - 1 - X \ge 0.$$ Then 
\begin{itemize}
\item $\xi$ ends in an odd number of 0s if $X$ does, 
\item $\xi$ ends in an even number of 0s, namely zero, if $X$ ends 
in a strictly positive even number of 0s. 
\end{itemize}
\end{lemma}

\begin{pf}
Recall that,
\begin{align}\label{fibrule}
 2F_1 = F_2 \text{ and } 2F_n = F_{n+1} + F_{n-2},
\end{align}
for $n\ge 2.$ 
We have that $Y := F_{2n} - 1 = (10)^{n - 1}1$. 
At first suppose that $X = 10^{2t+1}$, $t\in \N$. Then, if $t=0$, $\xi$ ends 
in three 0s, otherwise it ends in precisely one 0. 
Otherwise we must have, in Zeckendorf coding, 
$X = x010^{2t+1}$, $t\in \N$ for some bit-string $x > 0$. 
so we need to study an expression of the form 
$Y - X = (10)^{n-1} - x010^{2t+1}$. 
The trick we have in mind 
is probably easiest seen via an example: Put $n = 6$ so that 
$Y = F_{12} - 1 = 10101010101$ and 
suppose that $X = x0100000$,which gives $t=2$ . Then, repeated 
application of (\ref{fibrule}) give 
\begin{align*}
\xi &= 10101010101 - x0100000\\ 
&= 1201010101 - x0100000\\ 
&= 1112010101 - x0100000\\ 
&= 1111120101 - x0100000\\ 
&= z020101\\ 
&= z100201\\ 
&= z101010. 
\end{align*}
This number ends in an odd number of 0s irrespective of $z$, 
namely precisely one. This trick holds for all $n$ and $X$ which satisfy 
the conditions of the Lemma, except if $X$ ends in precisely one 0.
Hence we need to study this case separately. Here we get 
$$\xi = 1^{2(n-1)}2 - x010 = z102,$$ where $x$ and $z$ are bit-strings 
in Zeckendorf coding both with the least position at the 4th digit. 
Then, if $z$ ends in a 1, we get $\xi = w01110 = w10010$, for some $w$, 
which ends in precisely one 0. So assume that $z$ ends in a 0. 
Then, for some $w$, $\xi = w0102 = w1000$, so that, 
by Proposition \ref{prop:1}, 
$\xi$ ends in an odd number of 0s. In conclusion, the first item holds.

For the second item, notice that $\xi$ will end in zero 0s unless, 
in the subtraction, digit 2 gets a carry and (we may assume that $X$ is 
Zeckendorf coded) $X$ ends in 3 or more 0s. By assumption this 
has to be 4 or more 0s. But then, we do not need to add `a 
carry' to the second digit in the subtraction. 
\hfill \bqed
\end{pf}
\medskip

Notice that the second item does 
not hold if we exchange 'strictly positive' for 'non-negative'.
But, as will become apparent, it is only the first item 
of Lemma \ref{lemma:1} which is needed in 
the proof of Theorem \ref{theorem:2}. 

\begin{lemma}\label{lemma:2} Let $n\in \N$ and suppose that the integer 
$0 < X < F_{2n-1}$ ends in an even number of 0s, but not in $10^{2t+1}1$, 
$t\in \NN$. Then $$\varphi := F_{2n-1} - X$$ ends in an odd number of 0s.

Suppose that $T\in \N $ ends in $10^{2t + 1}1$ and $X$ in $10^{2s + 1}10$, 
where $T > X$ and $s, t\in \NN$. Then $$\xi := T - X$$ ends in an 
odd number of 0s.
\end{lemma}
\begin{pf}
For the first part, there are two cases to investigate, $X$ is either 
of the form 
\begin{enumerate}[(i)]
\item $x010^{2t}$, or 
\item $x010^{2t}1$,
\end{enumerate}
for some $t\in \N $.\\

\noindent Case (i): 
We have that $F_{2n - 1} = 10^{2(n-1)} = (10)^{s}110^{2t}$, where 
\begin{align}\label{nst}
n - 2 = s + t.
\end{align}
Then $$\varphi = F_{2n - 1} - X = (10)^s110^{2t} - x010^{2t} = r10^{2t+1},$$ 
(where $r = (10)^s0^{2(t+1)} - x0^{2(t+1)} > 0$) which, 
by Proposition \ref{prop:1} ends in an odd number of 0s independent of $r$.\\

\noindent Case (ii): We 
are going to prove that $\varphi$ ends in precisely one 0. 
By (\ref{nst}), $F_{2n-1} = (10)^{n-2}11$. Then 
$\varphi = (10)^{n-2}11 - x010^{2t}1 = r000(10)^t,$ which, by $t>0$, 
clearly ends in one 0, independent of $r$.\\ 

For the second part, as a first observation, if $T$ is of the form $y0101$, 
notice that, in Fibonacci numeration, $t\ge 0$ implies
\begin{align*}
\xi &= T - X\\ 
&= y0101 - x010^{2t+1}10\\ 
&= y0012 - x010^{2t+1}10\\ 
&= z002\\ 
&= z010, 
\end{align*}
for some bit-string $z$.
This idea generalizes to 
\begin{align*}
\xi &=  y010^{2s+1}1 - x010^{2t+1}10\\ 
&= y00(10)^s12 - x010^{2t+1}10\\ 
&= z010,
\end{align*} 
where, by assumption, $s\ge 0$.
Hence $\xi$ ends in precisely one 0 which resolves the second part of 
the lemma. \hfill \bqed
\end{pf}
\medskip

Suppose that the ordered pair $(X, Y)$ is of one of the forms 
in Theorem \ref{theorem:2} (ii) to (viii). Then 
Lemma \ref{lemma:1} and \ref{lemma:2} 
together imply that, if any of its option is $P$ then it has to be of 
the form 
\begin{align}\label{VW}
(V_i,W_i):=(X, Y)\ominus (B_i, A_i). 
\end{align}
(In other words, if $(X, Y)\ominus (A_i, B_i)$ is a legal option, it is $N$.)

But, by Theorem \ref{theorem:1} (iii), this is impossible if 
\begin{align*}
\frac{W_i}{V_i} > \phi .
\end{align*}

Hence, it suffices to investigate the cases 
\begin{align}\label{lephi}
\frac{W_i}{V_i} < \phi .
\end{align}
where, by Corollary \ref{cor:1}, both $W_i$ and $V_i$ end in 
an even number of 0s.

Before we prove Theorem \ref{theorem:2}, let us state a conjecture 
of `how far' we believe it could be extended by methods similar 
to those we have used in the Lemmas and below. 
(See also Figure \ref{figure:1} and Table \ref{table:1}.)

\begin{conj}\label{conj:1}
For all $n\ge 3$ and all $i$ such that 
\begin{itemize}
\item $A_i + B_i\le F_{2n-4}$, 
the position $(F_{2n-1}, F_{2n} - 1 - A_i - B_i)$ is $P$.
\item $A_i\le F_{2n-4}$, 
the position $(F_{2n-1} + A_i, F_{2n} - 1)$ is $P$.
\end{itemize}
\end{conj}

The following elementary result is an important tool for the proof 
of Theorem \ref{theorem:2}. 

\begin{lemma}\label{lemma:4}
For all $n\in \N$, 
\begin{align}\label{FFge}
\frac{F_{2n} - r}{F_{2n - 1} - s} > \phi  
\end{align}
if $0\le r \le \phi s$. 

For all $n\in \NN$, 
\begin{align}\label{FFle}
\frac{F_{2n+1} - r}{F_{2n} - s} < \phi
\end{align}
if $0\le \phi s \le r$.
\end{lemma}
\noindent{\bf Proof.} 
Notice that (\ref{FFge}) follows from, for all $n>0$, 
\begin{align}\label{FF1} 
\frac{F_{2n}}{F_{2n - 1}} > \phi  
\end{align}
 and (\ref{FFle}) from, for all $n\ge 0$, 
\begin{align}\label{FF2}
\frac{F_{2n+1}}{F_{2n}} < \phi,  
\end{align}
But (\ref{FF1}) and (\ref{FF2}) are easy. We give two alternative proofs. 
Clearly $\frac{F_1}{F_0} < \phi$. 
If $F_{2n+1} < \phi F_{2n}$ then $F_{2n+1} < \phi (F_{2n+2}-F_{2n+1})$ so that 
$F_{2n+1}\phi < F_{2n+2}$. If $F_{2n} > \phi F_{2n-1}$ then 
$F_{2n} > \phi (F_{2n+1}-F_{2n})$ so that $F_{2n}\phi > F_{2n+1}$. 

Another proof is given by using the well-known closed form expression 
$F_{n} = \frac{\phi^{n+1} - (1-\phi)^{n+1} }{\sqrt{5}}$. By this formula we get 
\begin{align*}
\frac{F_{n}}{F_{n-1}} &= \frac{\phi^{n+1} - (1-\phi)^{n+1} }{\phi^n - (1-\phi)^n }\\
&=  \phi\frac{\phi^{2n} - (-1)^{n+1} }{\phi^{2n} - (-1)^{n} },
\end{align*}
which gives the claim.
\hfill $\Box$\\

\noindent{\bf Proof of Theorem \ref{theorem:2}.}
For case (i) we need to prove that $(F_{2n} - 1, F_{2n} - 1)$ only 
has $N$-positions as options. By Lemma 
\ref{Fr82}, \ref{lemma:1} and \ref{lemma:2}, for all $i$, 
\begin{align}\label{F1b}
F_{2n} - 1 - b_i 
\end{align}
ends in an odd number of 0s and hence, by 
Proposition \ref{prop:1}, all options of the form 
$(F_{2n} - 1, F_{2n} - 1)\ominus (A_i, B_i)$ are $N$.\\

\noindent For case (ii), by (\ref{F1b}), we only need to be concerned 
with options of the form 
$(X, Y) := (F_{2n - 1}, F_{2n} - 1)\ominus (B_i, A_i)$. Notice that, by 
the first part of Lemma \ref{lemma:4}, for all 
$i\in \N$, $r = A_i + 1 \le B_i = s$, we get   
\begin{align}
\frac{Y}{X} = \frac{F_{2n}-r}{F_{2n-1}-s}
> \phi.
\end{align}
 Then Theorem \ref{theorem:1} (iii) gives the claim.\\

\noindent Case (iii). Here we want to prove 
that $(X, Y) := (F_{2n-1}, F_{2n} - 4)$ is $P$. 
The only move of the form $(b_i, a_i)$ which 
satisfies (\ref{FFle}) in Lemma \ref{lemma:4} is $(r,s) = (B_1, A_1) = (2, 1)$. 
By Theorem \ref{theorem:1} (iii) it then suffices to prove 
that $(X, Y)\ominus (2, 1)$ is $N$. 
In fact, we are going to demonstrate that 
$(X, Y)\ominus (2,1) = (F_{2n-1} - 2, F_{2n} - 5)$ 
has $(3, 3)$, which is $P$ (see also Table \ref{table:1}), as an option. 
For the latter it suffices to verify that
\begin{align}
M &:= (X,Y)\ominus (2,1) \ominus (100, 100)\notag\\ 
&= (10^{2t} - 10, 10^{2t+1} - 1000) \ominus (100, 100) \notag\\
&= (10^{2t} - 110, 10^{2t+1} - 1100)\label{plugin}
\end{align}
is a legal move (where $t = n-1$). 

Notice that, for $t\ge 2$, 
\begin{align*}
10^{2t} &= (10)^{t-3}100210\intertext{  and  } 10^{2t+1} &= (10)^{t-3}1002100.
\end{align*} 
By inserting these two identities into (\ref{plugin}) we get that $M$ is of the 
form $(z100, z1000)$ and hence legal.\\

\noindent Case (iv). By inspection, we have that 
$$\frac{F_{2n-1} - B_i}{ F_{2n} - 9 - A_i} > \phi ,$$ for all $i\ge 4$. 
So, by Lemma \ref{lemma:4}, it suffices to verify that each option 
$$(F_{2n-1}, F_{2n} - 9) \ominus \{(2, 1), (5, 3), (7, 4)\}$$ is $N$, 
respectively.
Hence, it suffices to demonstrate that
\begin{enumerate}[A:]
\item $(F_{2n-1} - 2, F_{2n} - 10)$ has the option $(6, 6)$, 
\item $(F_{2n-1} - 5, F_{2n} - 12)$ has the option $(F_{2n-2}, F_{2n-1} - 4)$, and
\item $(F_{2n-1} - 7, F_{2n} - 13)$ has the option $(3, 3)$,
\end{enumerate}
where the second item follows by case (ii). It suffices to 
verify that the moves are of the form in Lemma \ref{Fr82}. \\

Item A: We demonstrate that this is a legal move 
(using Fibonacci coding),
\begin{align*}
&(10^{2n}, (10)^{n-2}00101)\ominus (10, 1) \ominus (10100,10100)\\
&= (10^{2n}, (10)^{n-2}00101)\ominus (100000, 10101)\\
&= ((10)^{n-3}110000, (10)^{n-3}0110000)\ominus (100000, 10000)\\
&= ((10)^{n-3}010000, (10)^{n-3}0100000).
\end{align*}

Item B: Put $s = n-2$. Then, the move is
\begin{align*}
&(10^{2s+2} - 1000, 10^{2s+3} - 10101)\ominus (10^{2s+1}, 10^{2s+2} - 101)\\
&= (110^{2s} - 1000, 110^{2s+1} - 10000)\ominus (10^{2s+1}, 10^{2s+2})\\
&= (10^{2s} - 1000, 10^{2s+1} - 10000)\\
&= ((10)^{s-2}1100 - 1000, (10)^{s-2}11000 - 10000 )\\
&= ((10)^{s-2}0100, (10)^{s-2}01000),
\end{align*}
which is legal.\\

Item C: This is similar to item A:
\begin{align*}
&(10^{2n}, (10)^{n-2}00101)\ominus (1010, 101) \ominus (100,100)\\
&=((10)^{n-3}012011, (10)^{n-3}0101112)\ominus (10010, 1010)\\
&=((10)^{n-3}010100, (10)^{n-3}0101000).
\end{align*}

\noindent Case (v). By Lemma \ref{lemma:4} (as in Case (iii)), 
it suffices to demonstrate that 
$$(F_{2n-2} + 1, F_{2n-1} - 1)\ominus (2,1) = (F_{2n-2} - 1, F_{2n-1} - 2)$$ 
is of the form in Corollary \ref{cor:1}. 
But this hold since $F_{2n-2} - 1$ is of the form 
$10^{2t}-1=1010\ldots 1011-1=1010\ldots 1010$ which ends in precisely one 0.\\

\noindent Case (vi). By case (ii), positions of the form 
$(F_{2n-3},F_{2n-2}-1)$ is $P$. This also holds for $(11,11)$. 
Then one needs to verify that  
$$(F_{2n-1}+3, F_{2n}-1)\ominus (1,2)\ominus (11, 11)$$ 
and $$(F_{2n-1}+3, F_{2n}-1)\ominus (3,5)\ominus (F_{2n-3},F_{2n-2}-1)$$ are 
legal moves. We omit the details, since the methods are repetitions 
of the above. This suffice to prove the claim.\\

\noindent Case (vii). Let us demonstrate that the only options of 
$(F_{2n-1} + 4, F_{2n} - 1 )$ are of the form in Corollary \ref{cor:1}, that 
is, at least one of the coordinates ends in an odd number of 0s. By 
Lemma \ref{lemma:4}, we only need to check the moves 
$(1, 2)$, $(3, 5)$ and $(4, 7)$. In Fibonacci coding we have that 
$F_{2n}-1 = 1010\ldots 10101 = 1010\ldots01201 = 1010\ldots 01112 $. 
But then, by subtracting with $10, 1000$ and $1010$ respectively (and using 
the rule (\ref{fibrule})) we are done with this case.\\

\noindent Case (viii). By Lemma \ref{lemma:4}, here it suffices to verify 
that each one of the four options 
$(F_{2n-1} + 6, F_{2n} - 1)\ominus \{(2, 1), (5, 3), (7, 4), (10, 6)\}$ is $N$. 
We leave out much of the details since the verifications are repetitions 
of the above. However, a `rough line' goes as follows:

It may be verified that 
$(F_{2n-1} + 6, F_{2n} - 1)\ominus (1,2)$ 
has the $P$-position $(14, 14)$ as an option.
By Case (vi), $(F_{2n-3} + 3, F_{2n-2}-1)$ is $P$. This position is an option of 
$(F_{2n-1} + 6, F_{2n} - 1)\ominus (3, 5)$.
The option $(F_{2n-1} + 6, F_{2n} - 1) \ominus (4,7)$ is $N$. This follows by 
Corollary \ref{cor:1}, since $10^{2t}+10$ ends in an odd number of 0s.
Finally, it may be verified that 
$(F_{2n-1} + 6, F_{2n} - 1)\ominus (6,10)$ has the $P$-position 
$(1,1)$ as an option.
\hfill \bqed \\

The terminal positions of $\W^\star$ are all positions of the 
form $(0, n), n\in \N$. Denote the non-terminal $P$-positions of $\W^\star$
(with $a_i\le b_i$) by, in lexicographic order, $(a_1, b_1), (a_2, b_2),\ldots $
\begin{cor}\label{cor:2}
For $n\in \N$, 
$f(n) := \frac{b_n}{a_n}$ does not converge as $n\rightarrow \infty$. 
In particular, for all $\epsilon \in \R$ 
and $n\in \IZ_{>0}$ there is an $i \ge n$ such that
$f(i) = 1$ and a $j = j(\epsilon)\ge n$ such 
that $\phi-\epsilon < f(j) < \phi$.  
\end{cor}

\begin{pf}
This follows from Theorem \ref{theorem:2} (i) and (ii). In 
particular, notice that, for all $n\in \N$,
 $(F_{2n-1}, F_{2n}) $ is a $P$-position of Wythoff Nim.
\hfill \bqed
\end{pf}\medskip

\begin{ques}\label{quesWstar}
From Theorem \ref{theorem:2} it follows that we may characterize 
some $P$-positions of $\W^{\star}$ in polynomial time. Is there any method to 
extend these results to a polynomial time algorithm of determining if an 
arbitrary position is $P$?
\end{ques}

Numerical data, via computer simulations, motivate the following conjecture:

\begin{conj}\label{conj:2} Define the sets 
\begin{align*}
S_1 &:= \{3, 8, 11, 21, 32\},\\
S_2 &:= \{129, 362\},\\
S_3 &:= \{x\in \N\setminus\{19\}\mid \text{The Zeckendorf coding of 
$x$ ends in 101001}\},\\
S_4 &:= \{x\in \N\mid \text{The Zeckendorf coding of $x$ ends in 1}\}. 
\end{align*}
Then, the position $(i, i)$ belongs to $\P(\W^\star)$ if $i$ belongs 
to $(S_1\cup S_4)\setminus (S_2\cup S_3)$. 
It belongs to $\mathcal{N}(\W^\star)$ if $i$ belongs to 
$\N \setminus (S_1\cup S_2\cup S_3).$
\end{conj}

Let us give, in order of appearance, the Zeckendorf coding of the numbers 
in Conjecture \ref{conj:2} (they seem to have some special relevance 
to $\W^\star$ which I do not yet understand).
\begin{rem}
The Zeckendorf coding of
\begin{itemize}
\item $3, 8, 11, 21$ and $32$ are $100, 10000, 10100, 1000000$ 
and $1010100$ respectively,
\item $129$ and  $362$ are $1010001001$ and $101010001001$ respectively,
\item $19$ is $101001$.
\end{itemize}
\end{rem}

\section{Permutation games and the $\star$-operator}
We have defined the $\star$-operator in Section \ref{sec:1.1}.
Suppose that $(a_i)$ and $(b_i)$ are complementary sequences, both increasing.
Define $G$ by setting
$\M(G) := \{\{a_n, b_n\}\mid n\in \N \}$.
The Main Theorem in \cite{LHF} gives sufficient conditions on $a$ and $b$ 
(for example they may denote any pair of complementary Beatty sequences)
such that
\begin{align}\label{result}
\mathcal{P}(G^\star) = \M(G) \cup \{\boldsymbol 0\}
\end{align}
and therefore, 
\begin{align}\label{duality}
G^{\star\star} = G.
\end{align}
The question to try and classify all invariant subtraction games $G$ for 
which (\ref{result}) and (\ref{duality}) hold is left open. We will here 
ask a related, but more general question. Let us first explain what we mean 
by `convergence' of games.

For $n\in \NN$ and $G$ an invariant subtraction game, 
denote $$G_n = M(G)\cap\{(x,y)| x \le n\}.$$
Let $(G(k))$ denote a sequence of invariant subtraction games.

Suppose that there is a game $H$ such that, for all $n$, 
there is a $k$ such that $H_n = G_n(i)$ for all $i>k$.
Then $(G(k))$ converges and hence we can define the limit game 
$H = \lim_{k\in \N}G(k)$.

\begin{ques}\label{ques:1}
Let $G$ be an invariant subtraction game. Is it then true that the game 
$$H = \lim_{k\in\N} G^{2k}$$ exists?
\end{ques}

Clearly the answer to this question is affirmative for each game which 
satisfies (\ref{result}). But this is trivial, so we want to set out to try and 
find a larger family of games with an affirmative answer 
to Question \ref{ques:1}, but for which, in 
general, $G\ne G^{\star\star} $. We have defined our candidates, 
the permutation games, in Section \ref{sec:1.2}.

The first result is that the set of all permutation (involution) games 
is closed under the operation $\star\star$.

\begin{theorem}\label{thmbis}
If $G$ is a permutation game then, so is $G^{\star\star}.$
Furthermore, if $G$ is an involution game, so is $G^{\star\star}$.
\end{theorem}

Then we give and affirmative answer of Question \ref{ques:1} for $G$ a 
permutation game. 

\begin{theorem}\label{thmlim}
Let $G$ be a permutation game. Then $H = \lim_{k\in \N} G^{2k}$ exists. 
Furthermore, $H$ is a permutation game. If $G$ is an involution 
game, so is $H$. 
\end{theorem}

\begin{rem}
The involution games generalize the games defined 
by complementary sequences of integers discussed in \cite{LHF}. 
Suppose that $a = (a_i)$ and $b = (b_i)$ are 
two complementary sequences. Define 
$G$ by $\M(G) = \{\{a_i, b_i\}\}$. Then $G$ is an involution game. 
Let $G =$ 2-pile Nim. Then $G^\star$ is an involution game, 
but $\M(G^\star) = \{\{a_i, b_i\}\}$ gives $a_i = b_i$, for all $i$, so 
that $a$ and $b$ are not complementary. 
\end{rem}

\begin{rem}
The operator $\star\star$ may turn an invariant subtraction game 
which is not 
a permutation (involution) game into a permutation (involution) game. 
For example, define $G$ by 
$\M(G) = \{(i,i)\mid i\in \N\}\cup\{(1,2)\}$. 
Then, since $\P(G) = \M(\text{2-pile Nim})\cup\{\boldsymbol 0\}$, we get that 
$\M(G^{\star\star}) = \{(i, i)\mid i\in \N\}$.
\end{rem}
The next observation is proved in \cite[Lemma 2.2]{LHF}. Note here 
the importance of the word `invariant'.

\begin{lemma}[\cite{LHF}]\label{PM}
A move in an invariant 
subtraction game can never be a $P$-position\footnote{This requires 
normal play, see also Section \ref{sec:5.3}.}.
\end{lemma}

\begin{lemma}\label{terminal}
Suppose that $G$ is an invariant subtraction game on $\B(G) = \NN\times\NN$ 
which satisfies $\{\{0, x\}\mid x\in \NN \}\subset \P(G)$. Then 
\begin{enumerate}[(i)]
\item no two $P$-positions of $G^\star$ lie in the same row or column,
\item $\{\{0, x\}\mid x\in \NN\}\subset \P(G^{\star\star})$.
\end{enumerate}
\end{lemma}

\begin{pf}
Case (i) is obvious since $\M(G^\star) = \P(G)\backslash \{\boldsymbol 0\}$. 
For (ii), we apply the definition of $\star$ twice together with 
Lemma \ref{PM}.  Namely, the assumption 
$$\{\{0, x\}\mid x\in \N\}\subset \P(G)\backslash \{\boldsymbol 0\} 
= \M(G^\star)$$ implies
$$\{\{0, x\}\mid x\in \N\}\cap \P(G^\star ) = 
\{\{0, x\}\mid x\in \NN\}\cap \M(G^{\star\star}) = \emptyset .$$ 
Then each $(x, y)\in \M(G^{\star\star})$ satisfies $x > 0$ and $y > 0$. This 
gives that the set $\{\{0, x\}\mid x\in \NN\}$ is a subset of 
all terminal $P$-positions of $G^{\star\star}$, which gives (ii).
\hfill \bqed
\end{pf}
\medskip 

Since a permutation game satisfies the conditions of Lemma \ref{terminal}, 
we get the following corollary.

\begin{cor}\label{corterminal}
Let $G$ be a permutation game. Then $\{\{0, x\}\mid x\in \NN\}\subset \P(G)$ 
and $\{\{0, x\}\mid x\in \NN\}\subset \P(G^{\star\star})$.
\end{cor}
The next observation relaxes the requirements 
in Theorem \ref{theorem:1} (iii).
\begin{lemma}\label{finite}
Suppose that $G$ is an invariant subtraction game satisfying the assumptions 
in Lemma \ref{terminal}, for example a permutation game. Suppose further that 
the column $x$ contains the move $(x, y)$. Then $(x, z)$ is $N$ if $z > y$.
\end{lemma}
\noindent{\bf Proof.}
By the assumption, $z > y$ implies that $$(0,z - y) = (x, z)\ominus(x, y)$$ 
is $P$.
\hfill \bqed\\

\noindent{\bf Proof of Theorem \ref{thmbis}.}
Suppose that $G$ is a permutation game. We are going 
to show that the same holds for $G^{\star\star}$. First we demonstrate that, 
if there is a move in a column, then it is unique. Then we 
demonstrate that each column contains at least one move. The arguments 
obviously work fine with `column' exchanged for `row'. At last we prove 
the claim of `symmetry' in case $G$ is an involution game.\\

\noindent \emph{Uniqueness}: 
By the definition of a permutation game, no position 
of the form $\{0, x\}$ is a move in $G$. Hence all positions of this form 
belong to $\P(G) = \M(G^\star)\cup \{\boldsymbol 0\}$. But (by 
Lemma \ref{terminal} (i)) 
this gives that there can be no two positions of 
$\P(G^\star) = \M(G^{\star\star})\cup\{\boldsymbol 0\}$ in the same 
row or in the same column.\\

\noindent \emph{Existence}: 
Suppose that there is a least column, say $x_0 > 0$ 
which does not contain a move of $G^{\star\star}$. Then the 
set $$\{(x, y)\mid 0 < x < x_0\}\cap(\M(G^{\star\star})\cup \P(G^{\star\star}))$$ 
must be infinite. By the proof of `\emph{Uniqueness}' we already know 
that $$\M(G^{\star\star})\cap\{(x, y)\mid 0 < x < x_0\}$$ is finite. 
But all positions 
of the form $\{0, x\}$ are $P$-positions of $G$, hence, by 
Lemma \ref{terminal} (ii), also of $G^{\star\star}$. Then, by the minimality 
of $x_0$, $0 < r < x_0$, implies that column $r$ contains a move, 
say $(r, s)\in \M(G^{\star\star})$. Then, by Lemma \ref{finite}, 
$$\P(G^{\star\star})\cap\{(r, y)\mid 0 < r < x_0\}$$ is also finite. 
Hence there cannot exist such an $x_0$. By bijectivity 
of $\pi$, the argument works also with `columns' exchanged for `rows'.\\

\noindent \emph{Symmetry}: If $G$ is an involution game it follows that also 
$\M(G^\star) = \P(G)\setminus\{\boldsymbol 0\}$ is symmetric, which implies that 
$\M(G^{\star\star}) = \P(G^\star)\setminus\{\boldsymbol 0\} $ is symmetric. But, 
by the first part, we already know that $G^{\star\star}$ is a permutation game, 
hence also an involution game. 
\hfill \bqed
\medskip

\noindent{\bf Proof of Theorem \ref{thmlim}.}
Fix a permutation game $G$. Then, by induction on 
Theorem \ref{thmbis}, we get that, for all $k\in \NN$, $G^{2k}$ is also. 
Suppose now that, for a fixed $k\in \N$, there is a least column $x_0 > 0$ 
such that 
\begin{align}\label{mk}
(x_0, y) \in \M(G^{2k}) 
\end{align}
and 
\begin{align}\label{mk1}
(x_0, z) \in \M(G^{2k+2}), 
\end{align}
but $y\ne z$. (If there is no such $k$ then we are trivially done.) Then\\

\noindent \emph{Claim 1}: $(x_0, z)\in \mathcal{N}(G^{2k})$.\\

\noindent \emph{Claim 2}: For all $r < x_0$, $(r, s)\in \P (G^{2k})$ 
if and only if $(r, s)\in \P(G^{2k+2})$.\\

\noindent \emph{Claim 3}: $z > y$.\\

\noindent \emph{Claim 4}: $(x_0, z) \in \M(G^{2k+4})$. \\

\noindent {\bf Proof of Claim 1-4.} 
\noindent Claim 1: Suppose on the contrary that 
$(x_0, z)\in \mathcal{P}(G^{2k})$. 
Then $(x_0, z)\in \mathcal{M}(G^{2k+1})$, which (by \cite[Lemma 2.2]{LHF}) 
implies $(x_0, z)\in \mathcal{N}(G^{2k+1})$, which, 
by definition of $\star$, contradicts the assumption (\ref{mk1}).\\ 

\noindent Claim 2: Since, by assumption, the only move that differs 
in the two games is in column $x_0$ or greater, the claim follows.\\

\noindent Claim 3: By Claim 1, $(x_0, z)\in \mathcal{N}(G^{2k})$. Suppose  
that $z < y$. Then, by  (\ref{mk}) and the 
definition of a permutation game, there must exist a move, 
say $(r, s)\in \M(G^{2k})$, with $0 < r < x_0$ and $0 < s < z$ such that 
the option 
\begin{align}\label{XY}
(X, Y) &:= (x_0, z)\ominus (r, s)\notag\\
&\in\P(G^{2k}) = \M(G^{2k+1})\cup \{\boldsymbol 0\}.  
\end{align}
By minimality of $x_0$, we also have 
$$(r, s)\in \M(G^{2k+2}),$$ which, by definition of $\star$, implies 
$$(r, s)\in \P(G^{2k+1}).$$ This, together with (\ref{XY}), gives that 
\begin{align*}
(x_0, z) &= (X, Y)\oplus (r, s)\\
&\in\mathcal{N}(G^{2k+1}), 
\end{align*}
which, by definition of $\star$, contradicts the assumption (\ref{mk1}).
The claim follows.\\

\noindent Claim 4: Suppose this does not hold. 
Then, by Claim 3 and since, by Theorem 
\ref{thmbis}, $G^{2k+4}$ is a permutation game, there is a $w > z$ such that 
$(x_0, w)\in \M(G^{2k+4})$ and it is unique. By definition of $P$, 
it is clear that, in Claim 3, $z$ is the 
least number such that, for all legal moves 
\begin{align}\label{rs}
(r,s)\in \M(G^{2k+1}) 
\end{align}
we have that
$$(x_0, z)\ominus (r, s) \in \mathcal{N}(G^{2k+1}).$$ 

Observe that, by assuming that Claim 4 does not hold, we get 
that $(x_o,z)\in \n(G^{2k+3})$. 
Then, there has to be a move $(u, v)\in \M(G^{2k+3}) $ 
not of the form in (\ref{rs}), such that 
\begin{align}\label{x0zuv}
(x_0, z)\ominus (u, v)\in \P(G^{2k+3}).
\end{align}
But then 
\begin{align}\label{uv}
(u, v)\in \P(G^{2k+2})\setminus\{\boldsymbol 0\} 
\end{align}
so that Claim 2 together with (\ref{rs}) give 
$(u, v)\in \P(G^{2k})\setminus\{\boldsymbol 0\}$ which gives 
$(u, v)\in \M(G^{2k+1})$. But then (\ref{rs}) and the 
definition of $(u,v)$ gives a contradiction. 
Since, by Theorem \ref{thmbis}, $G^{2k+4}$ is a permutation game we are 
done with this case.

Again, by Theorem \ref{thmbis}, $G^{2k}$ is a permutation game 
for all $k\in \NN$. 
Together with Claim 4, this gives the existence of the permutation game $H$.

If $G$ is an involution game, then, by Theorem \ref{thmbis}, 
for all $k\in \NN$, $G^{2k}$ is. This gives that $H$ is also an involution game.
\hfill \bqed
\medskip
\section{Ornament games and complementary Beatty sequences with rational moduli}

In \cite{F2010} the `variant' game 'The Mouse game' is introduced. 
It answered the question: 
\begin{ques}\label{quesAviezri}
Is there an impartial game $G$ with the set of $P$-positions 
defined by complementary Beatty sequences with rational moduli? 
\end{ques}
The set of 
$P$-positions of the Mouse game is
\begin{align}\label{S}
S := \{\{a_i, b_i\}\mid i\in \N \}\cup \{\boldsymbol 0\}, 
\end{align}
where 
\begin{align}\label{Pmouse}
a_i := \left\lfloor\frac{3i}{2}\right\rfloor \text{ and } b_i := 3i - 1. 
\end{align}
The Mouse game is an extension of Wythoff Nim, where the available moves 
depend on which particular position of the board the next player moves from. 
(It is `variant', because it is not invariant.)
Precisely, if the position is $(x, y)$, where $y - x \equiv 0 \pmod 3$, 
then a player may move $$(x, y)\rightarrow (w, z),$$ 
provided $x - w \ge 0$, $y - z \ge 0$ and 
$$\mid (x - w) - (y - z) \mid \ \le 1.$$ 
Otherwise the moves are as in Wythoff Nim. 

In \cite{LHF} we gave an affirmative answer to the following question:
\begin{ques}\label{quesLHF}
Is there an invariant subtraction game $G$ with the set of (non-zero) 
$P$-positions defined by a pair of complementary Beatty sequences 
with rational moduli?  
\end{ques}
However we did not provide any closed form of the moves of any such game.
If, in Question \ref{quesLHF}, `rational' is exchanged for `irrational' 
the game of Wythoff Nim provided a solution over 100 years ago.

\subsection{The invariant moves of the Mouse trap}\label{trapsection}
Here we study an invariant subtraction game, 
\emph{the Mouse trap} = (the Mouse game)$^{\star\star}$, 
introduced in \cite{LHF}, with (by \cite[Main Theorem]{LHF}) an identical 
set of $P$-positions as the Mouse game. See Figure \ref{figure:3}. 
We present a closed form for the moves of this game and give an 
affirmative answer to this question.
\begin{ques}\label{quesPoly}
Is there an invariant subtraction game $G$ with the set of $P$-positions 
defined by complementary Beatty sequences with rational moduli 
and with the set of moves given by a polynomial time algorithm? 
\end{ques}
In fact, in Section \ref{subside} 
we give an affirmative answer to the same question, but this time 
for an infinite family of games. But the Mouse trap has an 
interest for its own sake. Here we give its invariant moves. 

\begin{figure}[!ht]
\centering
\includegraphics[width=0.7\textwidth]{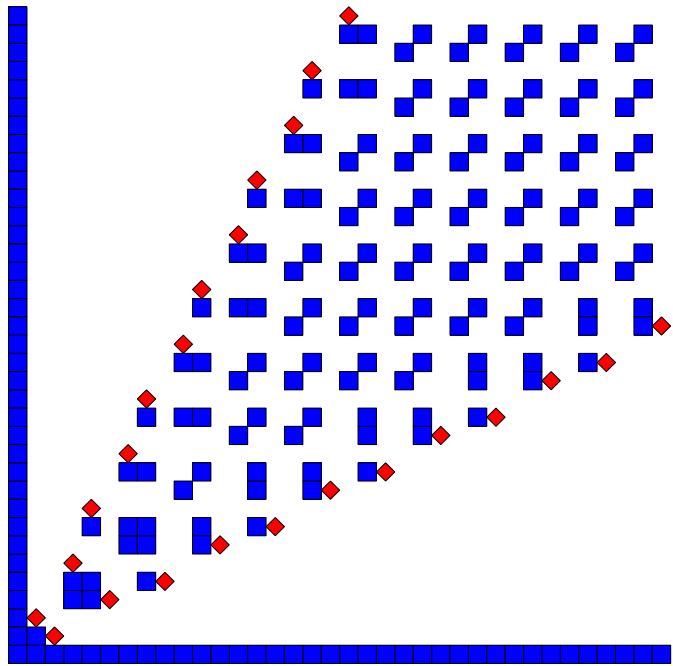}
\caption{The red diamonds and the 
blue squares represent the initial moves and $P$-positions of 
(Mouse game)$^\star$ respectively. Hence the blue squares, 
except $\boldsymbol 0$, symbolizes the moves of 'the Mouse trap' defined 
in Theorem \ref{theoremTrap}. The corresponding $P$-positions are then the 
red diamonds together with $\boldsymbol 0$. 
}\label{figure:2}
\end{figure}

\begin{theorem}\label{theoremTrap}
With notation as in (\ref{Pmouse}), define
\begin{align*}
\M_0 &= \{(1, 1), (3, 3), \{3, 4\}, (4, 4), \{4, 7\}, 
(6, 6), \{6, 7\}, \{6, 10\}\}\\
\M_1 &=  \{\{a_{2n-1}, a_{2m-1}\} \mid m,n\in \N, 3\le n\le m < 2n-1\}\\
\M_2 &= \{\{a_{2n}, a_{2m}\} \mid m,n\in \N, 3\le n\le m < 2n-2\}\\
\M_3 &= \{\{a_{2n},a_{4n-1}\}, \{a_{2n},a_{4n-3}\}\mid 3\le n \in \N \}\\
\M_4 &= \{\{0,x\}\mid x\in \N\}.
\end{align*}
Define $G =$'the Mouse trap' by 
\begin{align}\label{moveTrap}
\M(G) = \bigcup_{i=0}^4 \M_i.
\end{align}
Then, with notation as in (\ref{S}),
$$\P(G) = S.$$
\end{theorem}

\begin{pf}
Notice that, by (\ref{Pmouse}), for all $n\in \N$, 
\begin{align}
a_{2n}&= 3n,\label{3n}\\  
a_{2n-1}&= 3n - 2 \label{3n1}\intertext{ and }
b_{n}&= 3n - 1.\label{3n2}
\end{align}

It is easy to verify that, given the moves in $\M_0$, we get that 
$(a_1, b_1) = (1, 2), (a_2, b_2) = (3, 5), (a_3, b_3) = (4, 8)$ 
and $(a_4, b_4) = (6, 11)$ are the unique $P$-positions above 
the main diagonal up to and including column 6. (See also Figure \ref{figure:3}.)
So assume that the column is $\ge 7$. Let us begin with the direction \\

\noindent $N\rightarrow P$:
 Let $(X, Y)\not\in S$, 
with say $X\le Y$. Then we need to prove that $(X, Y)$ has an option in $S$. 
For three distinct classes of positions this is already clear, namely if 
\begin{itemize}
\item $X = b_i$, some $i$,
\item $X = a_i$ and $Y > b_i$, some $i$,
\item $(X, Y)\in \M(G)$.
\end{itemize}
The first two items follow immediately from the definition of $\M_4$
(see also \cite[Main theorem]{LHF}, Theorem \ref{theorem:1}). 
The third item is immediate by the definition of an invariant move, 
namely $(X, Y)\ominus (X, Y) = \boldsymbol 0$, which (by normal play) is $P$. 
(See also Lemma \ref{PM}.) So assume that $(X, Y)$ does not belong to any 
of these three classes of positions. Then, by the first two items, 
\begin{align}\label{XYac}
(X, Y) = (a_i, a_j) 
\end{align}
for some 
$i \le j$ and with $a_j < b_i$. But also, by the third item, $\M_1$ 
and $\M_2$, we have that 
\begin{align}\label{ij2}
i\not\equiv j \pmod 2, 
\end{align}
or equivalently
\begin{align}\label{aij3}
a_i\not\equiv a_j \pmod 3, 
\end{align}
\emph{except}, by definition of $\M_2$ and $\M_3$, if 
both $i$ and $j$ are even and $j\ge 2i - 4$, in which case 
$(X, Y)$ is of the form 
\begin{align}\label{third}
(a_{2n}, a_{4n-2})\text{ or }(a_{2n}, a_{4n-4}).
\end{align}

Altogether, we claim that:\\
 
\noindent Claim 1: If $X\equiv 0\pmod 3$ and $Y\equiv 1\pmod 3$ 
then the next player can move to the position $(5, 3)\in S.$\\ 

\noindent Claim 2: If $X\equiv 1\pmod 3$ and $Y\equiv 0\pmod 3$ 
then the next player can move to the position $(3, 5)\in S.$\\ 

\noindent Claim 3: If (\ref{third}) holds  
then the next player can move to the position $(3, 5)\in S.$\\

By Figure \ref{figure:3} it is not hard to justify the columns 7 to 11. 
For columns greater than 11, we begin with\\ 

\noindent Proof of Claim 1: We have that 
$(X, Y) = (3n+6, 3m+7)\ominus (5, 3) = (3n+1, 3m+4)\in \M_1$ if and only if 
$(X, Y)\prec (a_i, b_i)$, some $i$. We will demonstrate that the 
`worst possible case' gives $(X, Y)\ominus (5, 3) = (a_i, b_i - 1)$, some i. 
This happens whenever $Y/X$ is `maximized', 
that is, by definition of $\M_2$ and $\M_3$, whenever 
\begin{align*}
(X, Y) &= (a_{2n}, a_{4n-5})\\ 
&= (a_{2n}, a_{2(2n-2)-1})\\ 
&= (3n, 3(2n-2)-2)\\ 
&= (3n, 6n-8). 
\end{align*} 
Hence, we get 
\begin{align*}
(X, Y)\ominus(5, 3) &= (3n, 6n-8)\ominus (5, 3) \\
&= (3(n-1)-2, 6(n-1)-5) \\
&= (a_{2(n-1)-1}, b_{2(n-1)-1}-1), 
\end{align*}
as claimed. Here we have made repeated use of (\ref{3n}), (\ref{3n1}) 
and (\ref{3n2}).\\

\noindent Proof of Claim 2: $(3n+4, 3m+6)\ominus (3, 5) = (3n+1, 3m+1)\in \M_1$ 
if and only if $(X, Y)\prec (a_i, b_i)$, some $i$. As in Claim 1, the worst 
possible case gives $(X, Y)\ominus (5,3) = (a_i, b_i - 1)$. 
We omit the details.\\

\noindent Proof of Claim 3: $(3n+3, 3m+6)\ominus (3, 5) = (3n, 3m+1)\in \M_3$, 
We have either $(a_{2n}, a_{4n-2})\ominus (3, 5)$ or 
$(a_{2n}, a_{4n-4})\ominus (3, 5)$. The first case is equivalent to 
$(3n, 3(2n-1))\ominus (3, 5)= (3(n-1), 6(n-1)-2) 
= (a_{2(n-1)}, b_{2(n-1)}-1)\in \M_3$, since $a_{4n-1}=6n-2=b_{2n}-1$. The second 
case is treated in analogy to this.\\

\noindent $P\rightarrow N$: Assume $(X, Y)= (a_i,b_i)$, some $i$. It will 
become apparent that also for this case it suffices to analyze a 
`worst case scenario'. This is when the move is of the form $(a_j,b_j-1)$. 
We get the option 
\begin{align*}
(Z,W):=(X,Y)\ominus (a_j,b_j-1) = (a_i-a_j,b_i-b_j+1)=(a_i-a_j, 3(i-j)+1).
\end{align*} 
There are four cases to investigate, 
\begin{enumerate}[(i)]
\item $i=2n, j=2m,$,
\item $i=2n, j=2m-1,$
\item $i=2n-1, j=2m,$
\item $i=2n-1, j=2m-1$.
\end{enumerate}
where $ n, m\in \N$. For both cases (i) and (iv), we get the option 
\begin{align}
(Z, W) &= (3(n-m), 6(n-m)+1)\\ 
&= (a_{2(n-m)}, b_{2(n-m)}+2)\\
&\not\in S.
\end{align}
Case (ii) gives $3(n-m)+2=Z\equiv 2\pmod 3$, which is of the form of 
a `$b$-coordinate'. This gives $(Z,W)\not\in S$ 
since $Z < W = 6(n-m)+4$. For case (iii) we get that 
\begin{align*}
(Z,W) &= (3(n-m)-2, 6(n-m)-2)\\
&=(a_{2(n-m)-1},b_{2(n-m)-1}+2)\\
&\not\in S.
\end{align*}
By symmetry, we are done.
\hfill \bqed 
\end{pf} 
\medskip 

\begin{figure}[!ht]
\centering
\includegraphics[width=0.7\textwidth]{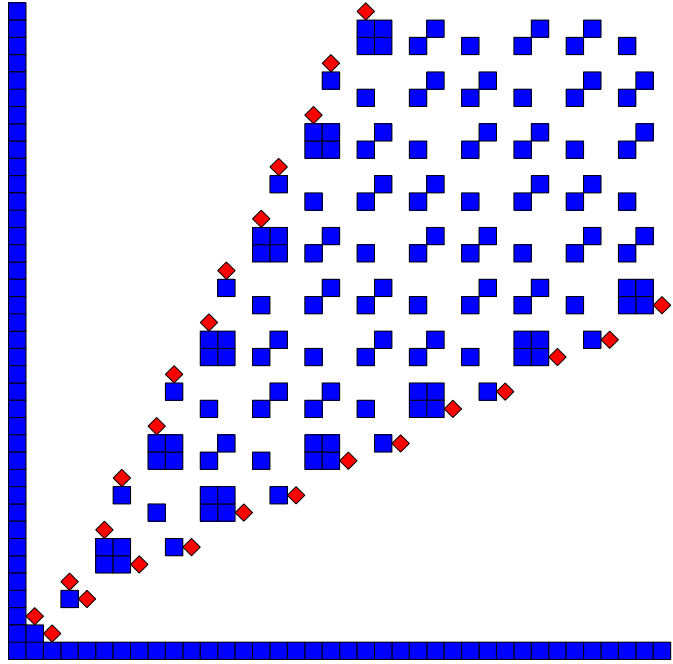}
\caption{An initial view of a 
close relative to (the Mouse game)$^\star$, namely where the 
complementary Beatty sequences are defined by 
$\alpha = 2/3, \delta =  -1/3, \beta = 1/3, \gamma = 2/3 $, 
notation as in (\ref{bea1}) and (\ref{bea2}). Notice that there are only 
three ornament games of Class 3. Thus the 
Figures \ref{figure:2}, \ref{figure:3} 
and \ref{figure:4} together illustrate this whole class. See also 
\cite{seminar} for updates of plots on more classes 
of ornament games.}\label{figure:3}
\end{figure}

\begin{rem}
Since the above $(a_i)$ and $(b_i)$ are complementary 
Beatty sequences, by \cite[Main Theorem]{LHF}, 
the moves defined in Theorem \ref{theoremTrap} 
are identical to the non-zero $P$-positions of (Mouse game)$^\star$. 
\end{rem}
\begin{rem}
By Theorem \ref{theoremTrap}, the decision problem for (Mouse game)$^\star$ 
has polynomial complexity. In contrast we do not know if this holds 
for (Wythoff Nim)$^\star$, where we, so far, only have the partial results 
in Corollary \ref{cor:1} and Theorem \ref{theorem:2} and together with the 
conjectures in Conjecture \ref{conj:1} and \ref{conj:2}.
\end{rem}
\subsection{Counting rational Beatty sequences and classes of ornament games}
In \cite{OB2003} the author provides a simple proof for the conditions of 
pairs of Beatty sequences to be complementary. To this purpose the sequences 
in (\ref{beatty1}) are translated to the forms 
\begin{align}\label{bea1}
\left(\left\lfloor \frac{n-\delta}{\alpha}\right\rfloor\right)_{n\in  \N} 
\end{align}
and 
\begin{align}\label{bea2}
\left(\left\lfloor \frac{n-\gamma}{\beta}\right\rfloor\right)_{n\in \N},
\end{align} 
$\alpha,\beta,\delta,\gamma\in \R$. 
Clearly, for density reasons, complementarity implies 
\begin{align}\label{alphabeta}
\alpha + \beta = 1 
\end{align}
and we may assume that $\alpha \le \beta$.
Hence, if one of the sequences has 
a rational modulus, then the other has also. 
Here we consider \emph{the rational case}. Then there is a 
least integer $q$ such that 
\begin{align}\label{qalpha}
q\alpha \in \N 
\end{align}
and by (\ref{alphabeta}), $q > 1$. Then, by \cite{OB2003}, 
the sequences are complementary if and only if 
\begin{align}\label{require1}
\frac{1}{q}\le \alpha + \delta \le 1 
\end{align}
and 
\begin{align}\label{require2}
\left\lceil q\delta \right\rceil + \left\lceil q\gamma \right\rceil = 1.
\end{align} 
Fix a constant $1 < C\in \N$ and let us estimate the number of pairs 
of rational 
Beatty sequences which together satisfy (\ref{alphabeta}), (\ref{qalpha}), 
(\ref{require1}), (\ref{require2}) and $q = C$. Let $\Xi$ denote 
the number of pairs of such sequences. It turns out that $\Xi$ only depends 
on $q$ (in particular it is independent of the reals $\delta$ and $\gamma$). 
\begin{prop}
Given $C\in \N$, $\Xi$ is finite. $\Xi = C\times (\varphi(C)+1)/2$, 
where $\varphi$ denotes the number 
of positive integers coprime with and less than $C$.
\end{prop}
\begin{pf}
With notation as above, by definition of $C=q$, there is a 
$p\in \{1, 2, \ldots , q-1\}$ such that 
$\alpha = \frac{p}{q}$ with $\gcd(p,q)=1$. Thus, we may 
rewrite (\ref{bea1}) as 
$(\lfloor \frac{qn - q\delta}{p}\rfloor)_{n\in  \N}$ (and (\ref{bea2}) as 
$(\lfloor \frac{qn - q\gamma}{q-p}\rfloor)_{n\in  \N}$).
For a fixed $\alpha$, put 
$$r_n(\delta) := \left\lfloor \frac{qn - q\delta}{p}\right\rfloor.$$ 

\noindent Claim 1: For 
$t\in \{-p+1,-p+1,\ldots , q-p\}$ and $s\in \R$, 
$r_n(\frac{t}{q}) < r_n(\frac{s}{q})$ implies that $s\le t-1$.\\

\noindent Claim 2: For each $n\in \N$, there exists an 
$m\in\{n, n + 1,\ldots , n + q - 1\}$ such 
that $r_m(\frac{t}{q}) < r_m(\frac{t-1}{q})$.\\

Assume that these two claims hold. Then 
we get that there are precisely $q$ distinct pairs of complementary
Beatty sequences satisfying (\ref{require1}) and  (\ref{require2}). 
Namely, we obviously need  
$\beta = 1 - \alpha$ and with $t$ as in Claim 1, 
take $\delta = \frac{t}{q}$ and $\gamma = \frac{1-t}{q}$.  
By symmetry (this is the division by 2), the proposition follows. 
But we need to prove the claims.\\ 

\noindent{Proof of Claim 1:} 
Suppose that 
$r_n(\frac{t}{q}) = \lfloor \frac{qn - t}{p}\rfloor 
< \lfloor \frac{qn - s}{p}\rfloor $. Then 
$\frac{qn - s}{p}-\frac{qn - t}{p} \ge \frac{1}{p}$ so that $-s+t \ge 1$.\\

\noindent{Proof of Claim 2:} We have that $\gcd(p,q) = 1$. Then there exists 
an $m\in \N$ such that $\frac{qm-t+1}{p}\in \N$. This implies that 

\begin{align*} 
\left\lfloor \frac{qm-t}{p}\right\rfloor &<\frac{qm-t}{p} + \frac{1}{p} \\
&= \left\lfloor \frac{qm-t+1}{p}\right\rfloor
\end{align*}
\hfill \bqed
\end{pf}
\medskip

Let $(\{\alpha n +\delta)$ and $(\beta n +\gamma )$ denote two complementary 
sequences with rational moduli. Suppose that the invariant subtraction 
game $G$ is defined by 
$\M(G)=\{\{\alpha n +\delta,\beta n +\gamma\} \mid n\in \N\}$. Then we call 
$G^\star$ an \emph{ornament game} (e.g. Figure \ref{figure:5}). If, 
in addition, $\alpha = \frac{C}{D}$ 
with $\gcd(C,D) = 1$, then $G^\star$ is an ornament game of \emph{Class} $C$.

\begin{cor}
The number of distinct ornament games of Class $C\in \N$ 
is $C\times (\varphi(C)+1)/2$. 
\end{cor}

\subsection{A subsided family of ornament games}\label{subside}
\begin{figure}[!ht]
\centering
\includegraphics[width=0.7\textwidth]{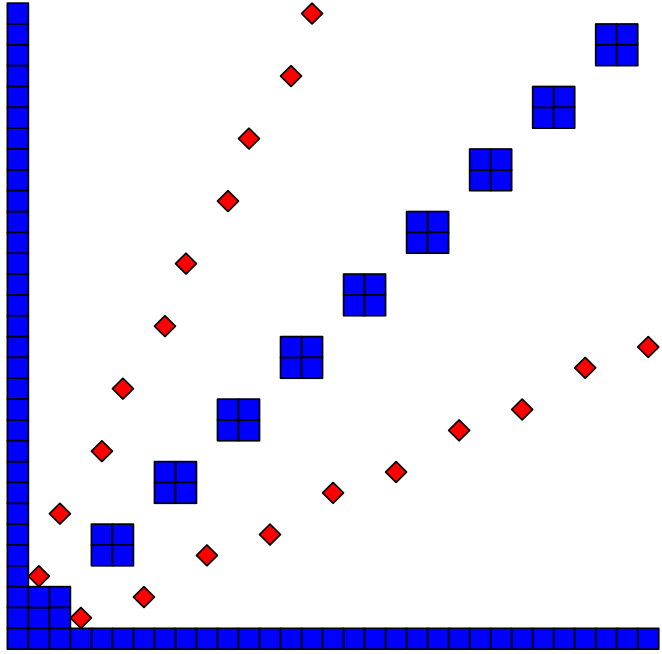}
\caption{An initial view of a 
close relative to (the Mouse game)$^\star$, namely where 
$\alpha = 2/3, \delta = 1/3, \beta = 1/3, \gamma = 0,$ see (\ref{bea1}) 
and (\ref{bea2}). 
The blue squares except $\boldsymbol 0$ may also be viewed 
as the first few moves of the unique Class 3 game in the 
family of `subsided' ornament games, 
see Theorem \ref{theoSub}.}\label{figure:4}
\end{figure}

In Section \ref{trapsection} we studied a special case of 
the Class 3 ornament games. 
With notation as in Theorem \ref{theoSub}, the games of the form 
$G^\star$ makes up a general family of ornament games with precisely 
one member in each class---thus, for example, the whole of Class 2 which 
consists of one single game $G^\star$ 
(with $\M(G) = \{\{2n, 2n-1\}\mid n\in \N\}$ and 
$\M(G^\star)\cup \boldsymbol 0 = 
\P(G) = \{\{2n-1, 2n-1\}\mid n\in \N\}\cup \{\{0,x\}\mid x\in \N\}$). 
We call $\{G^\star\}$ the family of 'subsided' ornament games 
(e.g. Figure \ref{figure:4}). 

\begin{theorem}\label{theoSub}
Let $2\le q\in \N$ and, for all $n\in \N$, put 
$$a_n := \left\lfloor\frac{qn-1}{q-1}\right\rfloor,  
b_n := qn.$$ Define $G$ by 
$$\M(G) = \left\{\left\{a_n, b_n\right\} \mid n\in \N\right\}.$$ 
Then $\P(G) = S\cup \{\{0,x\}\mid x\in \NN\}$, where 
\begin{align}\label{ornamentP}
S:=\{(qn+s, qn +t)\mid s,t\in \{1,2,\ldots , q-1\}, n\in \N \}.
\end{align}
\end{theorem}
\begin{pf}
Since, for all $n, a_n\le b_n$ and both sequences are increasing, by 
Theorem \ref{theorem:1} it suffices to study positions of the 
form $(a_i, a_j)$. By the same theorem and by symmetry, it suffices to 
study positions of the form $(x, y)$ with $a_n\le x\le y < b_n$. 
By elementary algebra we get that, 
for $i\in \{0,1,\ldots ,q-2\}$ and $n\in \N$, 
$$a_{(q-1)n-i} = qn-i-1$$ and $$b_n = qn.$$ 
Fix $r, s\in \N$ with $r\le s$, $i, j\in \{0, 1, \ldots , q-2\}$, and 
consider the position 
$$(X,Y):=\left(a_{(q-1)r - i}, a_{(q-1)s - j}\right) = 
\left(qr-i-1, qs-j-1\right).$$

\noindent $P\rightarrow N$: If $r = s$ then none of the options of $(X,Y)$ 
is of the form in (\ref{ornamentP}).\\ 
 
\noindent $N\rightarrow P$: If $r<s$ then there is an option of $(X,Y)$ of 
the form in (\ref{ornamentP}).\\

Together, these two claims suffices to prove the theorem. In an attempt 
to avoid uneccesary technicality we 
give the rest of the proof for the case $q = 3$. The general case may be 
treated in analogy\footnote{Incidentally, C.L.Bouton used the analog 
approach of 
`reducing technicality' in \cite{B1902} in proving the famous strategy of 
$q$-pile Nim. The case $q = 2$ is `too special' to be valid as a general proof 
and the cases $q\ge 3$ are analogous.}.\\

\noindent Proof of $P\rightarrow N$: For $q=3$ we get 
$$(X,Y) = (3r-i-1, 3r-j-1).$$ Since 
$b_1=3$, for the case $r=1$ we are done since $X$ has no legal options, so 
assume $r>1$. What is $(X,Y)\ominus(a_n, b_n)$? 
It suffices to prove that, for all $0<n\le r$, $X-a_n\ge 3(r-n)$. 
(By definition of $r$, all positions $(x, y)$ of the form in (\ref{ornamentP}) 
satisfy $x < 3(r-n)$, and $3n$ is the corresponding decrease of $Y$.) 
By definition of $a_n$, we have two cases to 
consider: $n = 2m$ or $n = 2m-1$, $m\in \N$. Suppose first that $n=2m$. 
Then, since $m>0$, $X - a_n = 3(r-m)-1 > 3r - 6m = 3(r-n)$. If, 
on the other hand $n = 2m - 1$, 
then $X - a_n = 3(r-m) \ge 3r-6m+3 = 3(r-n)$, with equality if and only if 
$m=1$. \\

\noindent Proof of $N\rightarrow P$: Here 
$$(X, Y) = \left(qr-i-1, qs-j-1\right),$$ with $r < s$ and $i, j\in \{0, 1\}$. 
It suffices to demonstrate the existance of 
a move of the form $(a_{2n-k}, b_{2n-k})$, $k\in \{0,1\}$ such that the option 
$$(Z,W):=(X,Y)\ominus (a_{2n-k},b_{2n-k})$$ is contained in the set $S$. By a 
simple calculation we get that 
\begin{align}\label{ZW}
(Z, W) = (3(r-n)-i+k, 3(s-2n+k)-j-1).
\end{align}
By the definition of the set $S$, it is required that $i\ne k$. Hence, 
we get two cases to consider:\\ 

\noindent Case $k < i$: (Clearly, for $q=3$ this forces $k=0$ and $i=1$, but 
in the coming we keep the symbols 
to make the  generalization to $q\ge 3$ more transparent.) 
Then both $-i+k$ and $-j-1$ are of the correct form as 
given in (\ref{ornamentP}). Thus, by (\ref{ZW}), it suffices to show that 
there is an $n\in \N$ satisfying $r-n = s-2n+k$. By the assumption $s>r$ and 
$k\in \{0,1\}$ we may take $n := s - r + k\in \N$.\\

\noindent Case $k > i$: For this case we may define $n$ and $k$ via: 
$$3(r-n) = 3(s-2(n-1)+k)+3$$ and $$0<k-i\le 3-j-1<3.$$ The former equation 
gives $n=s+k-r+1\in \N$ and the latter is clear. We are done.
\hfill \bqed
\end{pf}
\medskip

\begin{rem} Fix irrational moduli of a pair of complementary Beatty sequences. 
Then, by varying the real offsets, by \cite{OB2003}, 
it is easy to see that the corresponding (the irrational analogous 
to the ornament games) invariant subtraction games are uncountably many. 
Further classification of such games is left for future research. 
\end{rem}

 


%

\subsection{A new blocking maneuver on $k$-Wythoff Nim}
There is another, somewhat simpler, invariant game 
with the set of $P$-positions precisely $S$ as in (\ref{S}). 
For this game we define a certain 
'blocking maneuver/Muller twist' \cite{HL2006,L2010} 
on $k$-Wythoff Nim \cite{F1982}. 
The position $\boldsymbol 0$ does not appear among the 
pairs in (\ref{S}). Therefore, 
this position has been treated with some extra care in the following. 
Let $k\ge 3$. Define the game 'the constrained $k$-Mouse' as follows: 
Given a position $(x, y)\in \NN\times \NN$, move as in $k$-Wythoff 
Nim \cite{F1982}, that is a player may move 
$(x, y)\rightarrow (x - i, y - j)$, $0 \le \mid j - i \mid < k$, with 
$x\ge i$, $y\ge j$. But, before the next player moves, the previous player 
is allowed to block off at most 
$k - 2$ positions of the form $(x - i, y - j)$, 
\begin{align}\label{kmouse}
0 < \mid j - i \mid < k 
\end{align}
and declare that the next player may not move there. 
If the terminal position $\boldsymbol 0$ is of this form it may 
be blocked off, irrespective of the number of otherwise blocked off positions. 
When the next player has moved any blocked options are forgotten. 

\begin{theorem}
Let $G$ denote the constrained 3-Mouse. Then $\P(G) = S$. 
\end{theorem}
We omit the proof since it is in analogy to the blocking variations 
of Wythoff Nim presented in \cite{HL2006, L2010}. 
The move rules are considerably less technical for this 
blocking variation compared to the ones of 'the Mouse trap', both being 
invariant games. An open question is: What is $\P(G)$ for $k>3$?

\subsection{A conjecture on periodicity}
Let $S\subset \N \times \N$. Then $S$ is \emph{periodic} if, for 
all (sufficiently large) $(r, s)\in S$ implies that there is a pair
$(\alpha , \beta) \in \N\times \N$ such that, for all $n\in \NN $, 
\begin{align}\label{congr}
(r + \alpha n, s + \beta n)\in S.
\end{align} 
If, in a periodic set $S$, the number of distinct $(\alpha, \beta)$:s 
is bounded by a constant $k$, we say that $S$ is (at most) $k$-fold periodic. 

\begin{conj}\label{conj:4}
Let $(a_n) = (\lfloor \alpha n + \delta \rfloor)$ and 
$(b_i) = (\lfloor \beta n + \gamma \rfloor)$ denote complementary 
Beatty sequences, $\alpha ,\beta , \delta, \gamma \in \R$, $\alpha ,\beta >0$. 
Define $G$ via 
$\M(G) = \{(a_n, b_n) \mid n\in \N \}$, symmetric notation. Then $\P(G)$ 
is periodic if and only if the modulus $\alpha$ of $(a_n)$ is rational.
\end{conj}
Does this conjecture hold with `periodic' exchanged for `2-fold periodic'?\footnote{A related question is: Does each invariant 
subtraction game defined via complementary Beatty sequences with irrational 
moduli (such as $\W^\star$) have an infinite number of 
`$\log$-periodic' $P$-positions?}

\subsection{A question on Mis\`ere invariant subtraction games}\label{sec:5.3}
In \cite{LHF} (Remark 2) we have stated a belief that pairs of 
inhomogeneous complementary Beatty sequences might be worthy candidates 
as $P$-positions in Mis\`ere variations of invariant subtraction games. 
But having investigated a little further there seems to be some problem in 
constructing such games. The following example illustrates what can happen.
The sequence $$\{0,1\}, \{2,4\}, \{3,7\},  \{5,10\},\ldots $$ constitute 
the $P$-positions of the Mouse game minus $1$ in each coordinate. 
Hence the corresponding sequences of increasing integers are 
complementary Beatty sequences on the \emph{non-negative} integers. A 
Mis\`{e}re subtraction game requires that $\boldsymbol 0$ is $N$. Then 
$(0,1)$ is $P$ if it is a move. 
But the position $(2,3)$ has to be a move in a game for which the 
(complementary) $P$-positions begin with $(0,1), (2,4)$ and $(3,7)$. 
This follows since $(3,3)$ is $N$ and the only $P$-position $\prec (3,3)$ 
is $(0,1)$, which forces $(3,3)\ominus (0,1)$ to be a legal move.
Clearly, this `short-circuits' the $P$-positions $(2,4)$ and $(0,1)$. In 
this context it should be noted that in Lemma \ref{PM} which 
states that in an invariant subtraction game, a $P$-position can never 
be a move, normal play is required. Otherwise it is not true.
Rather, if $(X, Y)$ is a $P$-position in our above example, in Mis\`ere play,
the moves $$(X, Y)\ominus (0, 1)\text{ and } (X, Y)\ominus (1, 0)$$ 
are forbidden.

\begin{ques}
Are there Mis\`ere variations of invariant subtraction games 
on $\NN\times\NN$, such that the $P$-positions are given precisely by 
complementary inhomogeneous Beatty sequences of non-negative integers and with 
$(0,1)$ and $(1,0)$ \emph{the least} $P$-positions?
\end{ques}

\section{The dual of $k$-pile Nim}
Fix a $k\in \N$ and let $\mathcal{B} = \NN^k$. In this section we 
think of $\mathcal{B}$ as $k$ piles of tokens. Let $G$ denote $k$-pile Nim. 
Then $\M(G) = \{\{x, 0, 0, \ldots ,0\} \mid x \in \N\}$. 
It is well-known \cite{B1902} 
that $\P(G)$ is the set of all $k$-tuples with Nim-sum zero. 
We will now demonstrate that $k$-pile Nim has a dual game. 
\begin{theorem}\label{theorem:8}
Let $k\in \N$ and let $G$ denote $k$-pile Nim. Then $G^{\star\star} = G$, that 
is $G^\star$ is the dual game of Nim. 
\end{theorem}

\begin{pf}
Suppose that Alice plays first and Bob second in the game of $G^\star$, that is
the allowed moves are all non-zero $l$-tuples with Nim-sum zero, 
$2\le l\le k$. Notice 
that for this game a player must remove tokens from at least two piles.
It is then clear that the set of terminal positions is 
$\mathcal{T}(G^\star) = \M(G)\cup\{\boldsymbol 0\}\subseteq \P(G^\star)$. 
Denote the initial position with $(x_1, x_2, \ldots , x_k)$. We may assume 
$x_i\in \N$ for all $i$ and $x_1\le x_2\le \ldots \le x_k$. Then, a winning 
strategy for Alice is to, \emph{in one and the same move},  
\begin{itemize}
\item remove all $x_1$ tokens from pile $1$ and the same number 
from pile $x_2$,
\item if $x_2 > x_1$, remove $x_2 - x_1$ tokens from pile 2, that is 
the remaining ones, and the same number from pile 3, otherwise 
remove all tokens from pile 3 and $x_3$ tokens from pile 4,
\item continue in this manner until all piles are empty except possibly 
pile $k$ which now contains $x_k-x_{k-1}+\ldots +x_2-x_1\ge 0$ tokens if $k$ is 
even or $x_k-x_{k-1}+\ldots - x_2 + x_1\ge 0$ if $k$ is odd.
\end{itemize}
Alice's move is legal, since the Nim-sum of the number of removed tokens in 
the respective piles is zero. 
There is at most one pile of tokens left. Hence Bob cannot move and 
so Alice wins in her first move. This gives 
$\P(G^\star) = \mathcal{T}(G^\star)=\M(G)\cup \{\boldsymbol 0\}$. 
\hfill \bqed

\end{pf}

\begin{rem}
The proof of Theorem \ref{theorem:8} provides an intuitive winning strategy 
for \emph{$k$-pile Nim} without the mention of the concept 'Nim-sum'. 
Could, possibly, an averaged intelligent 
5 year old child learn to use this strategy? Shift the pile with the least 
number of tokens towards you, say an inch, and take the same number of tokens 
from any of the other piles and put them next to the pile you had 
just `shifted'. Continue in this same way with the remaining piles, 
until there is at most one pile of tokens left in the old place, `an 
inch further away'. The correct winning move, provided there is one, 
is to remove these tokens. We consider this strategy as only a little 
less intuitive than the rules themselves. At least it does not require any 
binary arithmetics.
\end{rem}

\noindent{\bf Acknowledgments.} Many thanks to Peter Hegarty, Aviezri Fraenkel 
and Richard Nowakowski for interesting discussions on topics related to 
the material in this paper.

\pagebreak
\begin{table}[ht]
\begin{center}
\begin{small}
\begin{tabular}{ l | c | c | c |c|| l | c |c| c|c | l | }
  $n$ & $a_n$ &fib.repr.& $b_n$ &fib.repr.& $n$ & $a_n$ &fib.repr.& $b_n$ &fib.repr. \\ \hline
  1   &1 &1&1  &1&  41 &33&1010101& 33 & 1010101\\\hline
  2   &3 &100&3  &100&42 &35&10000001& 35  &10000001\\\hline  
  3   &3 &100&4 &101 &43 & 37&10000100& 43 &10010001\\\hline
  4   &4 &101&4 &101 &44 & 38&10000101& 38 &10000101\\\hline
  5   &6 &1001&6 &1001 &45 & 40&10001001& 40 &10001001\\\hline
  6   &8 &10000&8 &10000 &46 & 40&10001001& 43 &10010001\\\hline
  7   &8 &10000&9 &10001 &47 & 42&10010000& 51 &10100101\\\hline
  8   &8 &10000&12 &10101 &48 & 42&10010000& 56 &100000001\\\hline
  9   &9 &10001&9  &10001&49 & 43&10010001& 43 &10010001\\\hline
  10  &9 &10001&12 &10101 &50 & 43&10010001& 46  &10010101\\\hline
  11  &11&10100& 11  &10100 &51 & 45&10010100& 48  &10100001\\\hline
  12  & 11&10100& 12 &10101 &52 & 45&10010100& 51 &10100101\\\hline  
  13  & 12&10101& 12 &10101 &53 & 45&10010100& 55  &100000000\\\hline 
  14  & 14&100001& 14 &100001 &54 & 46&10010101& 46 &10010101\\\hline
  15  & 16&100100& 17 &100101 &55 & 48&10100001& 48  &10100001\\\hline
  16  & 17&100101& 17 &100101 &56 & 48&10100001& 50  &10100100\\\hline
  17  & 19&101001& 19 &101001 &57 & 48&10100001& 51  &10100101\\\hline  
  18  & 19&101001& 21 &1000000 &58 & 48&10100001& 55 &100000000\\\hline
  19  & 21&1000000& 21 &1000000 &59 & 50&10100100& 59  &100000101\\\hline
  20  & 21&1000000& 25 &1000101 &60 & 51&10100101& 51  &10100101\\\hline
  21  & 21&1000000& 30 &1010001 &61 & 53&10101001& 55  &100000000\\\hline
  22  & 21&1000000& 33 &1010101 &62 & 53&10101001& 56  &100000001\\\hline  
  23  & 22&1000001& 22 & 1000001&63 & 53&10101001& 59  &100000101\\\hline
  24  & 22&1000001& 25 &1000101  &64 & 53&10101001& 63  &100010000\\\hline
  25  & 22&1000001& 33 &1010101 &65 & 55&100000000& 67 &100010101\\\hline
  26  & 24&1000100& 25 &1000101  &66 & 55&100000000& 72 &100100101\\\hline
  27  & 24&1000100& 27 &1001001 &67 & 55&100000000& 77 &101000001\\\hline
  28  & 24&1000100& 33 & 1010101&68 & 55&100000000& 80 &101000101\\\hline
  29  & 25&1000101 & 25 &1000101  &69 & 55&100000000& 85 &101010001\\\hline
  30  & 25&1000101 & 33 & 1010101&70 & 55&100000000& 88 &101010101\\\hline
  31  & 27&1001001& 27 &1001001 &71 & 56&100000001& 56 &100000001\\\hline
  32  & 27&1001001& 33 &1010101 &72 & 56&100000001& 59 &100000101\\\hline  
  33  & 29&1010000& 33 &1010101 &73 & 56&100000001& 67 &100010101\\\hline 
  34  & 29&1010000& 35 &10000001 &74 & 56&100000001& 88 &101010101\\\hline
  35  & 30&1010001& 30 &1010001 &75 & 58&100000100& 64 &100010001\\\hline
  36  & 30&1010001& 33 &1010101 &76 & 58&100000100& 67 &100010101\\\hline
  37  & 32&1010100& 32 & 1010100&77 & 58&100000100& 69 &100100001\\\hline  
  38  & 32&1010100& 33 &1010101 &78 & 58&100000100& 88 &101010101\\\hline
  39  & 32&1010100& 35 &10000001 &79 & 59&100000101& 59 &100000101\\\hline
  40  & 32&1010100& 42 &10010000 &80 & 59&100000101& 67 &100010101\\\hline
\end{tabular}
\end{small}
\end{center}\caption{The terminal $P$-positions of $W^\star$ are of the form 
$\{0, x\}, x\in \NN$. The first few $P$-positions 
which are not terminal are given here as $\{a_n, b_n\}$. We have included the
coding of these positions in the Fibonacci/Zeckendorf numeration system.}\label{table:1} 
\end{table}
\enddocument